\documentstyle[12pt,amssymb]{article}

\begin{document}

\newtheorem{thm}{Theorem}[section]
\newtheorem{lemma}[thm]{Lemma}
\newtheorem{defn}{Definition}[section]
\newtheorem{prop}[thm]{Proposition}
\newtheorem{corollary}[thm]{Corollary}
\newtheorem{remark}{Remark}

\renewcommand{\theequation}{\arabic{section}.\arabic{equation}}

\catcode`\@=11
\renewcommand{\section}{
        \setcounter{equation}{0}
        \@startsection {section}{1}{\z@}{-3.5ex plus -1ex minus
        -.2ex}{2.3ex plus .2ex}{\Large\bf}
}
\catcode`\@=12

\def\ee{\varepsilon}
\def\qed{\hfill\hbox{\vrule height1.5ex width.5em}}
\def\RR{{\bf R}}
\def\KK{{\bf K}}
\def\pf{\noindent{\bf Proof.} } 

\begin{doublespace}

\title{\Large \bf  Intrinsic Ultracontractivity,  Conditional Lifetimes
and Conditional Gauge for Symmetric 
Stable Processes on Rough Domains}
\author{Zhen-Qing Chen\thanks{The research
of this author is supported in part by NSA Grant
MDA904-98-1-0044}\\
Departments of Mathematics and Statistics\\ 
Cornell University \\
 Ithaca, NY 14853\\
Email: zchen@math.cornell.edu \smallskip \\
and  
\smallskip \\
Renming Song\\
Department of Mathematics\\
University of Illinois \\
Urbana, IL 61801\\
Email: rsong@math.uiuc.edu }
\date{ }
\maketitle

\begin{abstract}
For a  symmetric $\alpha$-stable process $X$ on $\RR^n$ with
$0<\alpha <2$, $n\geq 2$ and a domain $D \subset \RR^n$, let $L^D$ be the 
infinitesimal generator of the subprocess of $X$ killed upon 
leaving $D$. For a Kato class function $q$, it is shown
that $L^D+q$ is intrinsic ultracontractive on a 
H\"older domain $D$ of order 0. This is then used to
establish the conditional gauge theorem for
$X$ on
bounded Lipschitz domains in $\RR^n$. It is also shown that
the conditional lifetimes for symmetric stable process in
a H\"older domain of order 0 
are uniformly bounded.
\end{abstract}

\vspace{1in}


\vspace{.1in}

\noindent{\bf Running Title:} Conditional Gauge Theorem

\pagebreak

\section{Introduction}

A symmetric $\alpha$-stable process $X$ on $\RR^n$
is a L\'evy process
whose transition density $p(t, x-y)$
relative to Lebesgue measure is uniquely determined
by its Fourier transform 
$$
\int_{\RR^n} e^{i x\cdot \xi}  p(t, x) dx
=e^{-t |\xi|^\alpha} .
$$ 
Here $\alpha$ must be in the interval $(0, \, 2]$.
When $\alpha=2$, we get a Brownian motion running with a time clock
twice as fast as  the standard one.
Brownian motion  has been intensively studied
due to its  central role in modern probability theory
and its numerous important applications in other scientific areas
including many other branches of mathematics.
In the sequel, symmetric stable processes are referred to the case
when $0<\alpha<2$.
During the last thirty years, there has been an explosive
growth in the study of physical and economic systems that
can be successfully modeled with the use of stable processes.
Stable processes are now widely used in physics,
operation
research, queuing theory, mathematical finance and  risk estimation.
Recently some fine properties related to symmetric stable processes
and Riesz potential theory, as counterparts to Brownian motion and
Newtonian potential theory, have been studied, for example in
\cite{kr:bog}--\cite{s3:che}.

Let $X$ be a symmetric $\alpha$-stable process  
in $\RR^n$ with $n\geq 2$.
It is well known that $X$ is transient and has Green function
$G(x, y)= A(n, \alpha) |x-y|^{\alpha -n}$ where
\begin{equation}\label{eqn:2.1}
A(n, \alpha)=\frac{\alpha 2^{\alpha-1}\Gamma(\frac{\alpha+n}2)}
{\pi^{n/2}\Gamma(1-\frac{\alpha}2)}.
\end{equation}

\begin{defn}
A Borel measurable function $q$ on $\RR^n$ is said to be in the {\bf 
Kato class 
$\KK_{n, \alpha}$} if
\[
\lim_{r\downarrow 0}\sup_{x\in \RR^n}\int_{|x-y|\le r}
\frac {|q(y)|}{|x-y|^{n-\alpha}}dy=0.
\]
\end{defn}

For  $q \in \KK_{n,\alpha}$ and $t>0$, define
\[
e_q(t)=\exp\left( \int^t_0q(X_s)ds \right).
\]
For a domain $D \subset \RR^n$, the gauge function $g$ of $(D, q)$
is defined by $g(x)= E^x [ e_q(\tau_D )]$, $x \in D$ where
$\tau_D= \inf \{ t>0: \, X_t \notin D \}$. 
Chung and Rao \cite{r2:chu} showed the following {\it gauge theorem} 
for domain $D$ with 
finite Lebesgue measure:
the function $g$ is either identically
infinite or is bounded on $D$. In the latter case,
$(D, q)$ is said to be {\it gaugeable}
(with respect to $X$).  By Theorem 1 of Chung \cite{chung},
$(D, 0)$ is gaugeable for any domain $D$ with finite Lebesgue measure.
It is proved in \cite{s2:che}
that $(D, q)$ is gaugeable if and only if the first eigenvalue
of $L^D+q$ is negative, where $L^D$ is the non-positive definite
infinitesimal generator of the part process $X^D$ of $X$ killed
upon leaving domain $D$.

Now assume that $D$ is a bounded Lipschitz domain.
Recall that the Green function $G_D$ and Poisson kernel $K_D$ of $X$ in $D$
are determined by the following equations. For $x\in D$,
\begin{eqnarray*}
 E^x \left[\int_0^{\tau_D} f(X_s) ds \right] &=&\int_D G_D (x, y) f(y) dy
 \quad  \mbox{ for } f \geq 0 \mbox{ on } D , \\
 E^x \left[ \varphi (X_{\tau_D} ) \right] &=&  \int_{D^c} K_D (x, z)
    \varphi (z) dz  \quad \mbox{for } \varphi \geq 0 \mbox{ on } D^c.
\end{eqnarray*}
Fix an $x_0 \in D$, it is shown in \cite{s3:che} (see also
\cite{bog2}) that for each $z\in
\partial D$ and $x\in D$, 
\begin{equation}\label{eqn:martin}
 M_D(x, z)= \lim_{y\to z, \, y\in D}
{G_D(x, y)\over G_D(x_0, y)} \quad \hbox{exits and is finite}.
\end{equation}
$M_D$ is called the Martin kernel of $X$ in $D$. 

\begin{defn}\label{def:har}
Let $D$ be a domain in $\RR^n$.
A locally integrable function
$f$ defined on $D$ taking values in $(-\infty, \, \infty]$ and
satisfying the condition $\int_{\{|x|>1\}\cap D} |f(x)| |x|^{-(n+\alpha )}
dx <\infty$ is said to be 
{\bf superharmonic respect to $X^D$} if
$f$ is lower semicontinuous in $D$ and for each $x \in D$ and
each ball $B(x, r)$ with $\overline{B(x, r)}\subset D$,
$$
f(x)\ge E_x[f(X_{\tau_{B(x, r)}}); \tau_{B(x, r)}<\tau_D].
$$
A function $h$ is said to be {\bf harmonic  with respect to $X^D$}
if both $h$ and $-h$ are superharmonic respect to $X^D$.
\end{defn}

Suppose that $h\geq 0$ is superharmonic in $D$ with respect to $X^D$. 
Then by Theorem 2.3 of \cite{s3:che}, for
any domain $D_1\subset\overline{D_1}\subset D$,
$E_x [h^- (X^D_{\tau_{D_1}})] <\infty$ and
\[
h(x)\ge E_x[h(X^D_{\tau_{D_1}})] \quad \mbox{for every } x\in D_1. 
\]
This implies (cf. e.g. page 11 of Dynkin \cite{eb:dyn}) that
$$ h(x) \geq E^x [h(X^D_t)].
$$
Therefore one can define the $h$--conditioned stable process. 
Note that for each fixed $y\in D$ and
$z\in \partial D$, $G_D(x, y)$ and $M_D(x, z)$
are harmonic functions in $x$ with respect to  $X^{D\setminus \{y\}}$ and 
$X^D$,
respectively. For any $w\in\overline{D}^c$, $K_D(x, w)$
is superharmonic in $x$ with respect to $X^D$.
Therefore we can define the $G_D(\cdot, y)$--conditioned stable process,
$K_D(\cdot, w)$--conditioned stable process and the 
$M_D(\cdot, z)$--conditioned stable process.
The probability laws corresponding to these conditional stable processes
will be denoted by $P^x_y$, $P^x_w$ and $P^x_z$, respectively.
The following conditional gauge theorem is established in Chen and Song
\cite{s2:che}.

\begin{thm}\label{thm:cgauge}
Assume that $D$ is a bounded $C^{1,1}$ smooth domain in $\RR^n$,
$q\in \KK_{n, \alpha}$ and that $(D, q)$ is gaugeable.
Then exists a constant  $c>1$ such that
\[
c^{-1}\le\inf_{(x, z)\in D\times (\RR^n\setminus \partial D)}
E^x_z \left[e_q(\zeta )\right]
\le\sup_{(x, z)\in D\times (\RR^n\setminus \partial D) }
E^x_z \left[e_q(\zeta ) \right]
\le c.
\]
\end{thm}

Recently it is shown in \cite{s3:che} that under the condition
of Theorem \ref{thm:cgauge}, there is a constant $c>1$ such that
$$
c^{-1}\le\inf_{(x, z)\in D\times \partial D}
E^x_z \left[e_q(\zeta )\right]
\le\sup_{(x, z)\in D\times \partial D}
E^x_z \left[e_q(\zeta ) \right]
\le c.
$$

Gauge theorem for Brownian motion was  first proved by Chung and Rao in \cite{ra:chu} for bounded
$q$ and later was generalized to more general $q$
by various authors. The  conditional gauge theorem for Brownian motion
was first proved by Falkner \cite{ne:fal}
for bounded $q$ and a class of domains including bounded $C^2$ domains.
Extensions of this result to $q$ belonging to the Kato class and to bounded
$C^{1, 1}$ domains were given by Zhao in \cite{z1:zha} and \cite{z2:zha}.
The conditional gauge theorem has also been generalized
by Cranston, Fabes and Zhao \cite{fz:cra}
 to  diffusion processes whose infinitesimal generators are
uniformly elliptic divergence form operators
and to bounded Lipschitz domain in $\RR^n$.
For a more detailed story about gauge and conditional gauge theorem
for Brownian motion, the interested reader is referred to the recent book of Chung and Zhao \cite{zh:chu}.

\medskip

In this paper, we extend the above conditional gauge theorem
for symmetric stable processes
from bounded $C^{1,1}$ domains to bounded Lipschitz domains. 
We follow the idea from Chen and Song \cite{s2:che}, 
proving conditional gauge theorem through intrinsic ultracontractivity,
but with substantial improvements, motivated by 
Ba\~nuelos \cite{Ba}.
In \cite{s2:che}, sharp estimates on Green functions
of bounded $C^{1, 1}$ domains obtained in \cite{s1:che} 
were used to establish conditional gauge theorem. However these
estimates are no longer available for bounded Lipschitz
domains. We are able to circumvent it in this paper. 
We first show that, under the
assumption that $(D, q)$ is gaugeable, the Green
function $V_q$ of $L^D+q$ on $D$ is bounded by a constant multiple of
$G_D$. The latter is then used to prove the conditional gauge theorem
on bounded Lipschitz domains.
In order to establish $V_q\leq c \, G_D$ on $D\times D$ for some $c>0$,
we show that operator $L^D+q$ is intrinsic ultracontractive on $D$.
In fact we show that $L^D+q$ is intrinsic ultracontractive for any
H\"older domain $D$ of order $0$. For definitions of intrinsic ultracontractivity
and H\"older domain of order $0$,  see Definitions \ref{def:3.1} and
\ref{def:3.2} below.
We mention here that 
John domains, particularly bounded NTA domains and Lipschitz domains,
are  H\"{o}lder domains of order 0.
Under the assumption that $D$ has finite Lebesgue measure such that
$L^D$ is intrinsic ultracontractive, we show that there is a constant
$c>0$ such that for each non-trivial nonnegative superharmonic function
$h$ in $D$ with respect to $X^D$, $ \sup_{x\in D} E^x_h [ \tau_D ] \leq c $.
This especially implies that for a bounded Lipschitz domain $D$,
$$ 
   \sup_{x\in D, \, z\in \RR^n} E^x_z [\tau_D ] <\infty .
$$
Previously the boundedness of conditional lifetimes was proved for
bounded $C^{1,1}$ smooth domains for symmetric stable processes
in \cite{s1:che} and \cite{s3:che}. 
Conditional lifetime estimates for Brownian motion in planar domains
was first studied by M. Cranston and T. McConnell, in answering 
a question of K.~L. Chung. The first extension to several dimensions
for Brownian motion was done by M. Cranston, followed by many works
on important extensions to  more general domains and to 
elliptic diffusions (see
\cite{Ba}, \cite{BB} and the references therein).

\vspace{.1in}

\noindent{\bf Acknowledgement.}
We are grateful to Rodrigo Ba\~nuelos, Rich Bass and Chris Burdzy
for very helpful
discussions.
We  thank Krzysztof Bogdan
for very interesting and helpful discussions
at the MSRI, Berkeley, and University of Washington, Seattle.
Some of our results have been independently and simultaneously
obtained by him and will appear in \cite{bog3}.
Part of the research for this paper were conducted while the authors
were visiting MSRI in Berkeley. Financial support from MSRI are gratefully
acknowledged.

\section{Preliminaries}

Throughout this paper, we assume that $n\geq 2$ and 
$0<\alpha <2$.
Let $X$ be a symmetric $\alpha$-stable process  
in $\RR^n$.
It is well known that the Dirichlet form $({\cal E}, {\cal F}^{\RR^n})$
associated with $X$ is given by
\begin{eqnarray*}
{\cal E}(u, v)&=&A(n, \alpha)
\int_{\RR^n}\int_{\RR^n}\frac{(u(x)-u(y))(v(x)-v(y))}
{|x-y|^{n+\alpha}}\, dxdy\\
{\cal F}^{\RR^n}&=&
\left\{ u\in L^2(\RR^n):\int_{\RR^n}\int_{\RR^n}\frac{(u(x)-u(y))^2}
{|x-y|^{n+\alpha}}\, dxdy<\infty \right\},
\end{eqnarray*}
where $A(n, \alpha)$ is the constant in (\ref{eqn:2.1}).
As usual, we use $\{P_t\}_{t\geq 0}$
to denote the transition semigroup of $X$ and
$G$ to denote the potential of $\{P_t\}_{t\geq 0}$; that is, 
\[
Gf(x)=\int^{\infty}_0P_tf(x)dt.
\]

FROM now on, we assume that $D$ is a domain in $\RR^n$.
Adjoin an extra point $\partial$ to $D$ and set
\[
X^D_t(\omega)=\cases{ X_t(\omega) & if $t<\tau_D (\omega )$,\cr
\partial & if $t\ge \tau_D (\omega )$. \cr}
\]
The process $X^D$
 is called the symmetric $\alpha$-stable process killed upon leaving $D$,
or simply the killed symmetric $\alpha$-stable process on $D$. It is well
known (cf. \cite{ot:fuk})
 that the Dirichlet form corresponding to the killed symmetric
$\alpha$-stable process  $X^D$ on $D$ is $({\cal E}, {\cal F})$ where
\[
{\cal F}=\{u\in {\cal F}^{\RR^n}: 
\tilde{u}=0 \mbox{ quasi everywhere on } D^c\},
\]
where $\tilde{u}$ denotes a quasi continuous version of $u$.

We are going to use $P^D_t$ and $p^D(t, x, y)$
to denote the transition semigroup and transition density of $X^D$
respectively. $L^D$ will used to denote
the non--positive definite infinitesimal generator of $X^D$ on
$L^2(D, dx)$.
$G_D$, $K_D$ and $M_D$ will be used to denote the Green function,
Poisson kernel and Martin kernel
of $X$ on $D$ respectively. From \cite{s1:che} we know that
when $D$ satisfies the uniform exterior cone condition, $G_D$
and $K_D$ are related by
\begin{equation}\label{eqn:pois}
K_D(x, z)=A(n, \alpha)\int_D\frac{G_D(x, y)}{|y-z|^{n+\alpha}}dy,
\end{equation}
where $A(n, \alpha)$ is the constant in (\ref{eqn:2.1}).
Suppose that $h>0$ is a positive superharmonic function with respect to
$X^D$. We define
\[
p^D_h(t, x, y)=h(x)^{-1}p^D(t, x, y)h(y), \qquad t>0, x, y\in D .
\]
It is easy to check that $p^D_h$ is a transition
density and it determines a Markov process.
This process is called the 
$h$--conditioned symmetric stable process.

For any $y\in D$, the $G_D(\cdot, y)$--conditioned symmetric stable process
is a Markov process with state space $(D\setminus\{y\})\cup\{\partial\}$,
with lifetime $\zeta=\tau_{D\setminus\{y\}}$. We will use
$P^x_y$ and $E^x_y$ to denote the probability and expectation
with respect to this process.

For any $w\in \overline{D}^c$, the $K_D(\cdot, w)$--conditioned
symmetric stable process
is a Markov process with state space $D\cup\{\partial\}$,
with lifetime $\zeta=\tau_D$. We will use
$P^x_w$ and $E^x_w$ to denote the probability and expectation
with respect to this  process.

For any $z\in \partial D$, the $M_D(\cdot, z)$--conditioned
symmetric stable process
is a Markov process with state space $D\cup\{\partial\}$,
with lifetime $\zeta=\tau_D$. We will use
$P^x_z$ and $E^x_z$ to denote the probability and expectation
with respect to this  process.

In \cite{s1:che} we proved the following 3G Theorem.

\begin{thm}\label{T3g} 
Suppose that $D$ is a bounded $C^{1, 1}$ domain
in $\RR^n$ with  $n\ge 2$. Then
there exists a constant $c=c(D,n, \alpha)>0$ such that
\begin{eqnarray}
\frac{G_D(x, y)G_D(y, w)}{G_D(x, w)} &\le & \,
\frac{c\, |x-w|^{n-\alpha}}{|x-y|^{n-\alpha}|y-w|^{n-\alpha}}, 
  \quad x, y, w\in D,
\label{eqn:3g1}
\\
\frac{G_D(x, y)K_D(y, z)}{K_D(x, z)} &\le & \, 
 \frac{c\, |x-z|^{n-\alpha}}{|x-y|^{n-\alpha}|y-z|^{n-\alpha}}
,  \quad x, y\in D, z\in D^c.
\label{eqn:3g2}
\end{eqnarray}
\end{thm}

Using Theorem \ref{T3g} above and the scaling property we easily get the
following result.

\begin{corollary}\label{2a}
There exists a constant $c=c(n, \alpha)>0$ such that for any
ball $B$ in $\RR^n$ one has
\begin{eqnarray}
\frac{G_B(x, y)G_B(y, w)}{G_B(x, w)} &\le & \,
\frac{c\, |x-w|^{n-\alpha}}{|x-y|^{n-\alpha}|y-w|^{n-\alpha}}, 
  \quad x, y, w\in B,
\label{eqn:3g3}
\\
\frac{G_B(x, y)K_B(y, z)}{K_B(x, z)} &\le & \, 
 \frac{c\, |x-z|^{n-\alpha}}{|x-y|^{n-\alpha}|y-z|^{n-\alpha}}
,  \quad x, y\in B, z\in B^c.
\label{eqn:3g4}
\end{eqnarray}
\end{corollary}

3G Theorem actually holds for bounded Lipschitz domains.

\begin{thm}\label{2b}
Suppose that $D$ is a bounded Lipschitz domain in $\RR^n$. Then there
exists a constant $c=c(D, n, \alpha)>0$ such that
\begin{eqnarray}
\frac{G_D(x, y)G_D(y, w)}{G_D(x, w)}&\le&
\frac{c\, |x-w|^{n-\alpha}}{|x-y|^{n-\alpha}|y-w|^{n-\alpha}}, 
  \quad x, y, w\in D\label{eqn:3g5}\\
\frac{G_D(x, y)M_D(y, z)}{M_D(x, z)}&\le&
\frac{c\, |x-z|^{n-\alpha}}{|x-y|^{n-\alpha}|y-z|^{n-\alpha}}, 
  \quad x, y\in D, z\in \partial D.\label{eqn:3g6}
\end{eqnarray}
\end{thm}

\pf Note that the boundary Harnack inequality holds
on $D$ for positive harmonic functions in $D$ with respect to the
symmetric stable process $X$, due to Bogdan \cite{kr:bog}.
One can then prove (\ref{eqn:3g5}) by repeating the argument in Section
6.2 of \cite{zh:chu} so we omit the details here. 
Recalling (\ref{eqn:martin}), inequality (\ref{eqn:3g6}) follows from (\ref{eqn:3g5}) by passing $w\to z$.
\qed

We need the following result on the decomposition of Kato class
functions later on.

\begin{lemma}\label{kt:dcp}
Let $q$ have compact support. Then $q\in \KK_{n, \alpha}$
if and only if, for any $\ee >0$, there is a function
$q_\ee$ such that $q-q_\ee$ is bounded and 
\[
\sup_{x\in \RR^n}\int_{\RR^n}
\frac{|q_\ee (y)|}{|x-y|^{n-\alpha}}dy\le \ee
\]
\end{lemma}

\pf  The proof is the same as that of Theorem 4.16 in
\cite{sb:aiz}. \qed

\vspace{.2in}

In what follows, $q$ is an arbitrary but fixed
function in $\KK_{n,\alpha}$. 
For a domain $D$ in $\RR^n$, define 
\[
T_tf(x)=E^x\left[ e_q(t)f(X(t)); \, t<\tau_D \right], \qquad x\in D.
\]
The semigroup $T_t$ admits an integral kernel
$u_q(t, x, y)$ (cf. \cite{s2:che}).
The infinitesimal generator of the semigroup $T_t$ on $L^2(D, dx)$
is $L^D+q$. If $D$ has finite Lebesgue measure,  then
it is known (see Theorem 3.3 of \cite{s2:che}) 
that $L^D+q$ has discrete spectrum.
Let $\{\lambda_k, \,  k=0, 1, \cdots\}$ be
all the eigenvalues of $L^D+q$ written in decreasing
order, each repeated according to its multiplicity. Then
$\lambda_k\downarrow-\infty$ and the corresponding
eigenfunctions $\{\varphi_k, \,  k=0, 1, \cdots\}$ can be chosen that
they form an orthonormal basis of $L^2(D, dx)$. We know that
all the eigenfunctions $\varphi_k$ are bounded and the first
eigenfunction $\varphi_0$ can be chosen to be strictly positive in $D$
(cf. \cite{s2:che}).

\begin{defn}
A bounded Borel function $f$
defined on $\RR^n$ is said to a solution
to the equation 
\begin{equation}
(L^D+q)f(x)=0, \qquad x\in D\label{eqn:qd}
\end{equation}
if it is continuous in $D$ and for any open domain
$D_0\subset\overline{D_0}\subset D$,
\[
f(x)=E^x[e_q(\tau_{D_0})f(X_{\tau_{D_0}})], \qquad x\in D_0.
\]
\end{defn}

Clearly the first eigenfunction
$\varphi_0$ of $T_t$ is a positive solution of
\[
(L^D+q+\lambda_0)f(x)=0, \qquad x\in D.
\]
For positive solutions of (\ref{eqn:qd}) we have the
following uniform local Harnack inequality, which is 
applicable to $\varphi_0$ with $q+\lambda_0$ in place of $q$.

\begin{thm}\label{2c}
There exist two positive constants
$r_0=r_0(q)$ and $C=C(q)>0$ such that for any
solution $f$ of (\ref{eqn:qd}) which is strictly
positive on $D$ 
and for any ball $B(x_0, r)$ with
$0<r\le r_0$ and $B(x_0, 2r)\subset D$, one has
\[
\sup_{x\in B(x_0, r)}f(x)\le C\inf_{x\in B(x_0, r)}f(x).
\]
\end{thm}

\pf It follows from (\ref{eqn:3g4}) and the assumption of
$q\in \KK_{n, \alpha}$ that there exists a
positive number $R_0$ such that for any ball of radius $r\le R_0$
in $\RR^n$,
\[
\sup_{x\in B, z\in\overline{B}}\int_B 
\frac{G_B(x, y) |q(y)| K_B(y, z)}{K_B(x, z)}dy\le\frac12.
\]
By Jensen's inequality and Khas'minskii's lemma, 
\[
e^{-1/2}\le \inf_{x\in B, z\in\overline{B}^c}E^x_z[e_q(\tau_B)]
\le \sup_{x\in B, z\in\overline{B}^c}E^x_z[e_q(\tau_B)]\le 2.
\]
For any ball $B(x_0, r)$ with $0<r\le R_0/2$ and 
$\overline{B(x_0, 2r)}\subset D$, we have for any $x\in B(x_0, 2r)$,
\begin{eqnarray*}
f(x)&=&E^x[e_q(\tau_{B(x_0, 2r)})f(X_{\tau_{B(x_0, 2r)}})]\\
&=&\int_{B(x_0, 2r)^c}E^x_z[e_q(\tau_{B(x_0, 2r)})]\, f(z) \, 
K_{B(x_0, 2r)}(x,z) \, dz \\
&\leq & 2 E^x[ f(X_{\tau_{B(x_0, 2r)}})].
\end{eqnarray*}
By Harnack inequality
in \cite{cr:bas},  there exists a constant
$c=c(n, \alpha)>0$ such that
\[
\sup_{x\in B(x_0, r)}E^x[f(X_{\tau_{B(x_0, 2r)}})]\le c
\inf_{x\in B(x_0, r)}E^x[f(X_{\tau_{B(x_0, 2r)}})].
\]
Therefore for any $x, y\in B(x_0, y)$,
\begin{eqnarray*}
f(x)&\le &2c\, E^y[f(X_{\tau_{B(x_0, 2r)}})]\\
&\le&2c \, e^{1/2}
\int_{B(x_0, 2r)^c}E^y_z[e_q(\tau_{B(x_0, 2r)})] \, f(z) \, 
K_{B(x_0, 2r)}(y,z) \, dz\\
&=&2c\, e^{1/2}f(y),
\end{eqnarray*}
and the proof is now complete.
\qed

\vspace{.1in}

From the Theorem above we immediately get the following Harnack
inequality by a standard chain argument.

\begin{thm}\label{2d}
Suppose that $K$ is a compact subset of $D$. There exists a constant
$C=C(D, q, n, \alpha)>0$ such that for any solution
$f$ of (\ref{eqn:qd}) which is strictly
positive on $D$ one has
\[
\sup_{x\in K}f(x)\le C\inf_{x\in K}f(x).
\]
\end{thm}

\section{Intrinsic Ultracontractivity}

Let us first recall the definition of intrinsic ultracontractivity,
due to Davies and Simon \cite{si:dav}.
Suppose that $H$ is a semibounded self-adjoint operator on $L^2(D)$
with $D$ being a domain in $\RR^n$ and that $\{e^{Ht}, \, t> 0 \}$ is an irreducible
positivity-preserving semigroup with integral kernel $a(t, x, y)$. 
Assume that the top of the
spectrum $\mu_0$ of $H$ is an eigenvalue. In this case, $\mu_0$ has 
multiplicity one and the corresponding eigenfunction $\varphi_0$, 
normalized by $\|\varphi_0\|_2=1$, can be chosen to be positive almost everywhere on $D$.
$\varphi_0$ is called the ground state of $H$.

Let $U$ be the unitary operator $U$ from 
$L^2 \left(D, \varphi_0^2(x)dx\right)$ to
$L^2(D)$ given by $Uf=\varphi_0f$ and define $\widetilde H$ on 
$L^2 \left( D, \varphi_0^2(x)\,dx \right)$ by
\[
{\widetilde H}=U^{-1} \, (H-\mu_0) \, U.
\]
Then $e^{{\widetilde H}t}$ is an irreducible symmetric Markov semigroup
on $L^2\left(D, \varphi_0^2(x)\, dx \right)$ whose integral kernel with
respect to the measure $\varphi_0^2(x)dx$ is given by
\[
\frac{e^{-\mu_0 t}a(t, x, y)}{\varphi_0(x)\varphi_0(y)}.
\]

\begin{defn}\label{def:3.1}
 $H$ is said to be {\bf ultracontractive} if $e^{Ht}$ is a bounded
operator from $L^2(D)$ to $L^{\infty}(D)$ for all $t>0$.
$H$ is said to be {\bf  intrinsically ultracontractive} if  
${\widetilde H}$ is ultracontractive; that is, $e^{{\widetilde H}t}$
is a bounded operator from $L^2\left( D, \varphi_0^2(x)\,dx \right)$ to
$L^{\infty}\left( D, \varphi_0^2(x)\,dx \right)$ for all $t>0$.
\end{defn}

Ultracontractivity is connected to logarithmic Sobolev inequalities.
The connection between logarithmic Sobolev inequalities
and $L^p$ to $L^q$ bounds of semigroups was first given
by L. Gross \cite{gr1} in 1975.
E. Davies and B. Simons \cite{sb:aiz}
 adopted L. Gross's approach to
allow $q=\infty$ and therefore established the
connection between logarithmic Sobolev inequalities and
ultracontractivity.
(For an updated survey on the subject of logarithmic Sobolev inequalities
and contractivity properties of semigroups, see \cite{Bak}, \cite{gr2}.)
In \cite{Ba}, R. Ba\~nuelos proved the intrinsic ultracontractivity
for Shr\"odinger operators on uniformly H\"older domains of order
$\alpha \in (0, 2)$ using logarithmic Sobolev inequality characterization.
We will use the same strategy in this section; that is,
establishing intrinsic ultracontractivity
through logarithmic Sobolev inequalities.

In the rest of this section,  unless otherwise specified, 
$D$ is a domain in $\RR^n$ with finite Lebesgue measure, 
$q$ is a fixed
function in the Kato class $\KK_{n, \alpha}$.
 Recall that the semigroup $\{T_t, \, t>0 \})$ is
defined as follows.
\[
T_tf(x)=E^x\left[ e_q(t)f(X(t)); \, t<\tau_D \right], \qquad x\in D.
\]
$\{\lambda_k: k=0, 1, \cdots\}$ are
all the eigenvalues of $L^D+q$ written in decreasing
order, each repeated according to its multiplicity.
$\{\varphi_k: k=0, 1, \cdots\}$ are the corresponding
eigenfunctions, normalized so that they form an orthonormal
basis of $L^2(D, dx)$ and $\varphi_0>0$ on $D$.

The following result is proven in \cite{s2:che}.

\begin{thm}\label{ul:con}
The logarithmic Sobolev inequality holds for 
functions in $({\cal E}, {\cal F})$.
That is, for any $\eta>0$ and  $f\in {\cal F}
\cap L^{\infty}(D, dx)$, we have
\[
\int_Df^2\log|f|dx\le \eta{\cal E}(f, f)+\beta(\eta)\|f\|^2_2
+\|f\|^2_2\log\|f\|_2,
\]
with
\[
\beta(\eta)=-\frac{n}{2\alpha}\log\eta+c
\]
for some constant $c>0$.
\end{thm}

Recall that for any domain $D$ in $\RR^n$, the quasi--hyperbolic distance
between any two points $x_1$ and $x_2$ in $D$ is defined by
\[
\rho_D(x_1, x_2)=\inf_{\gamma}\int_{\gamma}\frac{ds}{\delta(x, \partial D)}
\]
where the infimum is taken over all rectifiable curves $\gamma$ joining
$x_1$ to $x_2$ in $D$ and $\delta(x, \partial D)$ is the Euclidean
distance between $x$ and $\partial D$.
Fix a point $x_0\in D$ which we call the center of $D$ and we may
assume without loss of generality that $\delta(x_0, \partial D)=1$.

\begin{defn}\label{def:3.2}
A domain $D$ in $\RR^n \, $ is a {\bf H\"older domain of order
0} if for a fixed $x_0\in D$, 
there exist constants $C_1$ and $C_2$ such that for all $x\in D$,
\[
\rho_D(x_0, x)\le C_1\log\left(\frac1{\delta(x, \partial D)}\right)
+C_2.
\]
\end{defn}
It is shown in Smith and  Stegenga \cite{st:smi} that
a  H\"older domain of order 0 is bounded.
It is well known that John domains, in particular 
bounded NTA domains and Lipschitz domains, are 
H\"{o}lder domains of order 0 (cf. \cite{Ba}).

\begin{lemma}\label{3a}
If $D$ is a H\"{o}lder domain of order 0, then for any $\beta>0$, there
exists a constant $C=C(D)>0$ such that
\[
\int_D(\rho_D(x_0, x))^{\beta}u^2(x)dx\le C{\cal E}(u, u), 
\qquad u\in {\cal F}
\]
\end{lemma}

\pf From \cite{st:smi} we know that for any $\beta>0$ we have
\[
\int_D(\rho_D(x_0, x))^{\beta}dx<\infty.
\]
It follows from the Sobolev inequality (see formula (1.5.20) of 
\cite{ot:fuk} or Theorem 1 on page 119 of \cite{el:ste}) that there
is constant $C_1>0$ such that for any $u\in {\cal F}$,
\[
\|u\|_{p_0}\le C_1\sqrt{{\cal E}(u, u)},
\]
where $p_0$ is such that $1/p_0=1/2-\alpha/(2n)$.
Let $p=n/(n-\alpha)$ and $p'=n/\alpha$. Applying H\"{o}lder's
inequality we get that for any $ u\in {\cal F}$,
\begin{eqnarray*}
\int_D(\rho_D(x_0, x))^{\beta}u^2(x)dx&\le&
\left(\int_D(\rho_D(x_0, x))^{\beta p'}dx\right)^{1/p'}\left(
\int_D |u|^{p_0}(x)dx\right)^{(n-\alpha)/n}\\
&\le&C_1\left(\int_D(\rho_D(x_0, x))^{\beta p'}dx\right)^{1/p'}
{\cal E}(u, u).
\end{eqnarray*}
The proof is now complete.
\qed

\begin{thm}\label{3b}
If $D$ is a  H\"{o}lder domain of order 0, 
then for any $\varepsilon>0$
and any $\sigma>0$ we have
\[
\int_Df^2\log\frac1{\varphi_0}dx\le\varepsilon
{\cal E}(f, f)+\beta(\varepsilon)\|f\|^2_2,
\qquad f\in {\cal F}
\]
with 
\[
\beta(\varepsilon)=c_1\varepsilon^{-\sigma}+c_2
\]
for some positive constants $c_1$ and $c_2$.
Here $\varphi_0$ is the ground state of $L^D+q$.
\end{thm}

\pf
Let $W=\{Q_j\}$ be a Whitney decomposition of $D$. This is a decomposition
of $D$ into closed cubes $Q$ with the following
three properties (see \cite{el:ste} for details)
\begin{description}
\item{(1)} for $j\neq k$, the interior of $Q_j$ and the interior of $Q_k$
are disjoint;
\item{(2)} if $Q_j$ and $Q_k$ intersects, then
\[
\frac14\le\frac{\mbox{diam}(Q_j)}{\mbox{diam}(Q_k)}\le 4;
\]
\item{(3)} for any $j$,
\[
1\le\frac{\delta(Q_j, \partial D )}{\mbox{diam}(Q_j)}\le 4.
\]
\end{description}
Let $x_0$ be a fixed point in $D$ and $x_0\in Q_0$. If $Q_k\in W$,
we say that $Q_0=Q(0)\rightarrow Q(1)\rightarrow\cdots\rightarrow Q(m)=
Q_k$ is a Whitney chain connecting $Q_0$ and $Q_k$ of length $m$ if
$Q(i)\in W$ for all $i$ and $Q(i)$ and $Q(i+1)$ have touching edges for
all $i$. We define the Whitney distance $d(Q_0, Q_k)$ to be the length
of the shortest Whitney chain connecting $Q_0$ and $Q_k$. If $x\in Q_k$
we define $\tilde{\rho}(x_0, x)=d(Q_0, Q_k)$. It is well known and easy
to prove that this distance is comparable with $\rho_D$, the quasi--hyperbolic
distance.

By Theorem \ref{2c}, property (3) of the Whitney
decomposition, the boundedness of $D$ and  the equivalence of
$\tilde{\rho}$ and $\rho_D$,  there is a constant $C_1=
C_1(D)>0$ such that for any $Q\in W$ we have
\[
\sup_{x\in Q}\varphi_0(x)\le C_1\inf_{x\in Q}\varphi_0(x).
\]
Therefore there exists a constant $C_2=C_2(D)>0$ such that
\begin{equation}
\varphi_0(x)\ge e^{-C_2\rho_D(x_0, x)}\varphi_0(x_0), \qquad x\in D.\label{eqn:eig}
\end{equation}
For any $p>1$, let $p'$ be its conjugate. By (\ref{eqn:eig}) and Lemma \ref{3a},
we get that for any $\varepsilon>0$ and $u\in {\cal F}$,
\begin{eqnarray*}
\int_Du^2\log\frac1{\varphi_0}dx&\le&C_2\int_D\rho_D(x_0, x)u(x)^2dx\\
&=&C_2\int_D\frac{\varepsilon^{1/p}}{\varepsilon^{1/p}}\rho_D(x_0, x)u(x)^2dx\\
&\le&\varepsilon C_2\int_D(\rho_D(x_0, x))^pu(x)^2dx+
C_2\varepsilon^{-p'/p}\int_Du^2(x)dx\\
&\le& \varepsilon \, C_3 \, {\cal E}(u, u)
+C_2 \, \varepsilon^{-p'/p}\int_Du^2(x)dx,
\end{eqnarray*}
where $C_3$ is a positive constant depending on $D$ only. The proof
is now complete.
\qed

Combining the two theorems above we get the
following result.

\begin{thm}\label{3c}
If $D$ is a H\"{o}lder domain of order 0, 
then for any $\varepsilon>0$
and any $\sigma>0$ we have
\[
\int_Df^2\log\frac{|f|}{\varphi_0}dx\le \eta{\cal E}(f, f)+\beta(\eta)\|f\|^2_2
+\|f\|^2_2\log\|f\|_2, \qquad f\in {\cal F}\cap L^{\infty}(D, dx)
\]
with 
\[
\beta(\varepsilon)=-\frac{n}{2\alpha}\log\varepsilon+
c_1\varepsilon^{-\sigma}+c_2
\]
for some positive constants $c_1$ and $c_2$.
\end{thm}

With the result above, we can easily get our main result
of this section.

\begin{thm}\label{iu}
Assume that $D$ is a H\"{o}lder domain of order 0.
Then
$L^D+q$ is intrinsically ultracontractive. More precisely,
\[
\frac{e^{-\lambda_0 t}u_q(t, x, y)}{\varphi_0(x)\varphi_0(y)}
\le e^{2M(t/2)} \quad \hbox{for all } x, y\in D \hbox{ and } t>0,
\]
where
\[
M(t)=\frac1t\int^t_0 A(\varepsilon)d\varepsilon
\]
with
\[
A(\varepsilon)=\cases{-\frac{n}{2\alpha}
\log\varepsilon+c_1\varepsilon^{-1/3}+c_2 &  for  $\varepsilon\le 1$, \cr
c_1+c_2 & for  $\varepsilon>1$. \cr}
\]
for some positive constants $c_1$ and $c_2$.
\end{thm}

\pf By taking $\sigma=1/3$ in Theorem \ref{3c},
we get that for any $\varepsilon>0$
and any $f\in {\cal F}\cap L^{\infty}(D, dx)$,
\begin{equation}
\int_Df^2\log\frac{|f|}{\varphi_0}dx\le \varepsilon{\cal E}(f, f)+
\beta_1(\varepsilon)\|f\|^2_2
+\|f\|^2_2\log\|f\|_2, \label{eqn:lsi}
\end{equation}
with
\[
\beta_1(\varepsilon)=-\frac{n}{2\alpha}\log\varepsilon+
c_1\varepsilon^{-1/3}+c_2
\]
for some positive constants $c_1$ and $c_2$.

Suppose that $(\widetilde{\cal{ T}},
\widetilde{\cal{ F}})$ is the
Dirichlet form on $L^2(m)$ with $m(dx)={\varphi}^2_0dx$ 
associated with the semigroup whose integral kernel with
respect to the measure $m$ is given  by 
\[
\frac{e^{-\lambda_0 t}u_q(t, x, y)}{{\varphi}_0(x){\varphi}_0(y)}.
\]
Then
\[
\widetilde{\cal F}=\{f: f{\varphi}_0\in {\cal F}\}
\]
and
\[
\widetilde{\cal T}(f, h)={\cal E}(f{\varphi}_0, h{\varphi}_0)
-\int_Dqf{\varphi}_0\, h{\varphi}_0\, dx+\lambda_0\int_Dfh\, dm
\]
Since $q\in \KK_{n, \alpha}$, by Theorem 3.2 of \cite{r2:son}
there exists a constant $B>0$ such that
\[
\int_D|q|u^2 dx\le\frac12 {\cal E}(u, u)+B\int_Du^2 dx, \qquad u\in {\cal F}.
\]
Thus
\[
\widetilde{\cal T}(h, h)\ge \frac12 {\cal E}(h{\varphi}_0, 
h{\varphi}_0)
-(B-\lambda_0)\int_Dh^2 dm.
\]
By putting $f=h{\varphi}_0$ in (\ref{eqn:lsi}) we get that for 
$h\in \widetilde {\cal F}\cap L^{\infty}(D, dm)$,
\begin{equation}
\int_Dh^2\log|h| \, dm\leq
 2\varepsilon \, \widetilde{\cal T}(h, h)+ \left(\beta_1(\varepsilon)
+2(B-\lambda_0)\right)
\int_Dh^2 dm+\|h\|^2_{L^2(m)}\log\|h\|_{L^2(m)}.
\label{s4:log}
\end{equation}
Therefore  for $0<\varepsilon\le 1$ and $h \in \widetilde {\cal F}
\cap L^{\infty}(D, dm)$,
\begin{equation}
\int_Dh^2\log|h| \, dm\le \varepsilon \, \widetilde{\cal{ T}}(h, h)+
\beta_2(\varepsilon)\int_Dh^2 dm
+\|h\|^2_{L^2(m)}\log\|h\|_{L^2(m)},
\label{s8:log}
\end{equation}
where
\[
\beta_2 (\varepsilon) =-{n\over 2 \alpha}\log \varepsilon +c_3 \varepsilon^{-1/3}+c_4,
\]
for some constants $c_3,\,  c_4>0$.
For $\varepsilon>1$,  since $\widetilde {\cal T}$ is nonnegative and 
(\ref{s8:log})
holds for $\varepsilon=1$, we have for any $h\in \widetilde {\cal F}
\cap L^{\infty}(D, dm)$,
\begin{eqnarray}
\int_Dh^2\log|h|dm&\le& \widetilde{\cal{ T}}(h, h)+\beta_1(1)\int_Dh^2 dm
+\|h\|^2_{L^2(m)}\log\|h\|_{L^2(m)}\nonumber \\
&\le&\varepsilon\widetilde{\cal{ T}}(h, h)+\beta_1(1)\int_Dh^2 dm
+\|h\|^2_{L^2(m)}\log\|h\|_{L^2(m)}.
\label{s9:log}
\end{eqnarray}
Combining (\ref{s8:log}) and (\ref{s9:log}) we get that for any $\varepsilon>0$
and any $h \in \widetilde {\cal F}
\cap L^{\infty}(D, dm)$,
\begin{equation}
\int_Dh^2\log|h|dm\le \varepsilon\widetilde{\cal{ T}}(h, h)+
A(\varepsilon)\|h\|^2_{L^2(m)}
+\|h\|^2_{L^2(m)}\log\|h\|_{L^2(m)},
\label{s10:log}
\end{equation}
with
\[
A(\varepsilon)=\cases{-\frac{n}{2\alpha}\log\varepsilon+c_5\varepsilon^{-1/3}+c_6 & 
for  $\eta\le 1$, \cr
c_5 +c_6 & for $\varepsilon>1$. \cr}
\]
for some positive constants $c_5$ and $c_6$.
By Corollary 2.2.8 of \cite{eb:dav}, we have
\[
\frac{e^{-\lambda_0 t}u_q(t, x, y)}{\varphi_0(x)\varphi_0(y)}
\le e^{2M(t/2)} \quad \hbox{for all } x, y\in D \, \hbox{ and } t>0,
\]
where
\[
M(t)=\frac1t\int^t_0 A(\varepsilon)d\varepsilon <\infty.
\]
\qed

Using the same argument as that for
Theorem 6 in R. Smits \cite{ro:smi}, we have

\begin{thm}\label{3d}
Assume that $D$ is a domain in $\RR^n$ with finite Lebesgue measure
such that $L^D+q$ is
intrinsic ultracontractive. Then
there exists $C>0$ such that for any $t>1$,
\[
e^{(\lambda_1-\lambda_0)t} \le \sup_{x, y\in D}
\left|\frac{e^{-\lambda_0 t}u_q(t, x, y)}{\varphi_0(x)\varphi_0(y)}-1 
\right|\le C
e^{(\lambda_1-\lambda_0)t}.
\]
\end{thm}

\section{Conditional Lifetimes}

Assume in this section that
 $D$ is a domain in $\RR^n$ with finite Lebesgue measure
such that $L^D$ is intrinsic ultracontractive, unless otherwise
specified. In particular, we know from
Theorem \ref{iu} that a H\"{o}lder domain of order 0
satisfies this assumption.
Since $(D, 0)$ is gaugeable, the first eigenvalue $\mu_0$ of $L^D$
is negative. Let
$\phi_0$ be the ground state of $L^D$. Recall that
$P_D$ is the transition density function for the killed
symmetric stable process $X^D$.
Similar to Corollary 1 of Ba\~nuelos \cite{Ba}, we have

\begin{thm}\label{thm:lifetime}
Under the assumption given at the beginning of this section.
There is a constant $c>0$ such that
\begin{description}
\item{(1)} $ e^{\mu_0 t} \phi_0 (x) \phi_0 (y) \leq p^D(t, x, y)
 \leq c e^{\mu_0 t} \phi_0 (x) \phi_o(y)$ for all $x, y \in D$ and $t>0$;
\item{(2)} Let ${\rm SH}^+$ denote all non-trivial nonnegative
superharmonic functions in $D$ with respect to
$X^D$. Then 
$$ \sup_{x\in D, \, h\in {\rm SH}^+} E^x_h [\tau_D ] <\infty ; $$
\item{(3)} For $h\in {\rm SH}^+$, 
$$ \lim_{t\to \infty } e^{-\mu_0 t} P^x_h (\tau_D >t)
= {\phi_0(x) \over h(x)} \int_D \phi_0(y) h(y) .
$$
In particular, $\lim_{t\to \infty } {1\over t} \log P^x_h (\tau_D >t)
=\mu_0$.
\end{description}
\end{thm}

\pf (1) It follows directly from Theorem \ref{3d}.

\noindent (2) Note that for each $h \in {\rm SH}^+$, by (1)
$$ h(x) \geq \int_D p^D(1, x, y) h(y) dy \geq e^{\mu_0} \phi_0 (x)
 \int_D \phi_0 (y) h(y) dy \quad \hbox{ for } x \in D.$$
Therefore
$$ \sup_{x\in D, \, h\in {\rm SH}^+} {\phi_0(x) \over h(x)}
\int_D \phi_0(y) h(y) \, dy \leq  e^{-\mu_0} <\infty .
$$
Therefore by (1)
\begin{eqnarray*}
&& \sup_{x\in D, \, h\in {\rm SH}^+} E^x_h [\tau_D ] \\
&=& \sup_{x\in D, \, h\in {\rm SH}^+}
{1\over h(x) }\int_0^\infty \int_D p^D(t, x, y) h(y) \, dy \\
&\leq & c \int_0^\infty e^{\mu_0 t} \, dt \,
\sup_{x\in D, \, h\in {\rm SH}^+} 
{\phi_0(x) \over h(x)} \int_D \phi_0(y) h(y) \, dy <\infty .
\end{eqnarray*}

\noindent (3) By Theorem \ref{3d}
\begin{eqnarray*}
&&\lim_{t\to \infty } e^{-\mu_0 t} P^x_h (\tau_D >t) \\
&=& \lim_{t\to \infty } e^{-\mu_0 t} \phi_0(x)^{-1}
   \int_D p^D (t, x, y) h(y) dy \\
&=& {\phi_0(x) \over h(x)} \int_D \phi_0(y) h(y) .
\end{eqnarray*}
\qed

When $D$ is a bounded Lipschitz domain, $G_D(x, y)$, $M_D(x, z)$ and
$K_D(x, w)$ are superharmonic functions in $x$ with respect to
$X^D$ for each fixed $y\in D$, $z\in \partial D$ and $w\in {\overline D}^c$,
respectively. The above theorem in particular implies that

\begin{corollary} (Conditional Lifetimes) 
Assume that $D$ is a bounded Lipschitz domain. Then
$$ \sup_{x\in D, \, z\in \RR^n} E^x_z [ \tau_D ] <\infty . $$
\end{corollary}

\section{Conditional Gauge Theorem}

Throughout this section,  $D$ is a bounded Lipschitz domain. Recall that 
$L^D$ is  the non-positive definite infinitesimal generator of the killed
$\alpha$-stable process on $D$. For $q\in \KK_{n, \alpha}$, let
$u_q(t, x, y)$ be the kernel of the following Feynman-Kac semigroup
\[
T_tf(x)=E^x\left[e_q(t)f(X_t)1_{\{t<\tau_D\}}\right].
\]
Note that the semigroup $T_t$ only depends on function $q$ through $q1_D$
so we may assume that $q=0$ off $D$.
The following result is proven in \cite{s2:che}.

\begin{thm}\label{sq:gat}
Suppose that $q\in \KK_{n, \alpha}$ is such that
\[
\sup_{x, y\in D}\int_D\frac{G_D(x, z)|q(z)|G_D(z, y)}{G_D(x, y)}dz\le
\frac12.
\]
Then we have
\begin{equation}
e^{-1/2}G_D(x, y)\le V_q(x, y)\le 2 G_D(x, y),\label{gc:smq}
\end{equation}
where
\[
V_q(x, y)=\int^{\infty}_0u_q(t, x, y)dt.
\]
\end{thm}

\begin{thm}\label{g1:cmp} 
Assume that $q\in \KK_{n, \alpha}$ and $(D, q)$ is gaugeable. Then
 there is a constant  $c>0$ such that
\[
V_q(x, y)\le c\,G_D(x, y) \quad \hbox{ for all } x, y\in D.
\]
\end{thm}

\pf 
By Lemma \ref{kt:dcp} and Theorem \ref{2b}, the function $q$ can be 
decomposed  as 
$q=q_1+q_2$ with $q_1$ bounded and 
$q_2\in \KK_{n, \alpha}$ satisfying
\begin{equation}\label{eqn:kato}
\sup_{x, y\in D}\int_D\frac{G_D(x, z)|q_2(z)|G_D(z, y)}{G_D(x, y)}dz\le
\frac12.
\end{equation}
Therefore by Theorem \ref{sq:gat} 
\begin{equation}
e^{-1/2}G_D(x, y)\le V_{q_2}(x, y)\le 2 G_D(x, y).\label{eqn:gcp}
\end{equation}

Let $\{\nu_k, \,  k=0, 1, \cdots\}$ be all the eigenvalues of $L^D+q_2$
written in decreasing order, each repeated according to its multiplicity.
and let $\{\psi_k, \,  k=0, 1,\cdots\}$ be the corresponding
eigenfunctions with $\psi>0$ on $D$. We can assume that $\{\psi_k, \,  k=0, 1,\cdots\}$
form an orthonormal basis of $L^2(D, dx)$.
Since by  Khas'minskii's lemma $E^x_y [ e_{q_2} (\zeta)]< 2$ for
$x, y \in D$,   we know from Theorem 3.11 of \cite{s2:che} that $\nu_0<0$. 

By Theorem \ref{3d}, there is a $t_1>1$ such that
for all $t\ge t_1$ and all $x, y\in D$,
\[
\frac12\le\frac{e^{-\nu_0 t}u_{q_2}(t, x, y)}{\psi_0(x)\psi_0(y)}
\le\frac32.
\]
Therefore for all $x, y\in D$,
\begin{eqnarray}
V_{q_2}(x, y)&=&\int^{\infty}_0u_{q_2}(t, x, y)dt
\nonumber\\
&\ge&\int^{\infty}_{t_1}u_{q_2}(t, x, y)dt
\nonumber\\
&\ge&\frac12\psi_0(x)\psi_0(y)\int^{\infty}_{t_1}e^{\nu_0 t}\, dt
\nonumber\\
&\ge& C_1 \, \psi_0(x)\psi_0(y),\label{eqn:gbd}
\end{eqnarray}
for some positive constant $C_1>0$.

Since $L^D+q=(L^D+q_2)+q_1$ and $q_1$ is bounded, it follows from 
Theorem 3.4 of \cite{si:dav} that the first eigenfunction
$\varphi_0$ of $L^D+q$ is comparable to the first eigenfunction
$\psi_0$ of $L^D+q_2$, i.\/e., there exists a constant $C_2>1$ such that
\begin{equation}
C^{-1}_2\psi_0\le \varphi_0\le C_2\psi_0.\label{eqn:ecp}
\end{equation}
By Theorem \ref{3d} again,  there is a $t_2>1$ such that
for all $t\ge t_2$ and all $x, y\in D$,
\[
\frac12\le\frac{e^{-\lambda_0 t}u_q(t, x, y)}{\varphi_0(x)\varphi_0(y)}
\le\frac32.
\]
Since the gauge function of $(D, q)$ is assumed to be finite,
the first eigenvalue $\lambda_0$ is negative by Theorem 3.11
of \cite{s2:che}. Therefore it follows from 
(\ref{eqn:ecp}), (\ref{eqn:gbd}) and (\ref{eqn:gcp})
 that for any $x, y\in D$,
\begin{eqnarray*}
\int^{\infty}_{t_2}u_q(t, x, y)dt&\le&\frac32\varphi_0(x)
\varphi_0(y)\int^{\infty}_{t_2}e^{\lambda_0 t}dt\\
&\le& C_3\varphi_0(x)\varphi_0(y)
\le C_4\psi_0(x)\psi_0(y)\\
&\le& C_5V_{q_2}(x, y)
\le C_6G_D(x, y),
\end{eqnarray*}
where $C_3, C_4, C_5, C_6$ are positive constants.
Since $q_1$ is bounded, 
\begin{eqnarray*}
\int^{t_2}_0u_q(t, x, y)dt&\le&e^{\|q_1\|_{\infty}t_2}\int^{t_2}_0
u_{q_2}(t, x, y)dt\\
&\le& e^{\|q_1\|_{\infty}t_2}V_{q_2}(x, y)\le 2e^{\|q_1\|_{\infty}t_2}
G_D(x, y).
\end{eqnarray*}
Hence there exists a constant $C>1$ such that
\[
V_q(x, y)\le CG_D(x, y), \qquad x, y\in D.
\]
\qed

Applying the Theorem above and the 3G Theorem we can get
the following results which correspond to Theorems 5.3--5.7
of \cite{s2:che}. The proofs are exactly the same as those
of the corresponding results in \cite{s2:che} so the details
are omitted.

\begin{thm}\label{gn:rel}
Assume that $(D, q)$ is gaugeable. Then for
 all $x, y\in D$ with $x\neq y$,
\begin{eqnarray}
V_q(x, y)&=&G_D(x, y)+\int_DV_q(x, u)q(u)G_D(u, y)du \label{g1:rel} \\
V_q(x, y)&=&G_D(x, y)+\int_DG_D(x, u)q(u)V_q(u, y)du.\label{g2:rel}
\end{eqnarray}
\end{thm}

\begin{thm}\label{ga:rep}
Assume that $(D, q)$ is gaugeable. Then for
 all $(x, y)\in D\times D$ with  $x\neq y$,
\[
 E^x_y \left[ e_q(\zeta )\right]
=1+G_D(x, y)^{-1}\int_DV_q(x, w)q(w)G_D(w, y)dw.
\]
\end{thm}

\begin{thm}\label{g2:rep}
Assume that $(D, q)$ is gaugeable. Then for
all $(x, y)\in D\times D$ with  $x\neq y$,
\[
E^x_y \left[ e_q(\zeta )\right] =\frac{V_q(x, y)}{G_D(x, y)}.
\]
\end{thm}

\begin{thm}\label{co:gau}
Assume that $(D, q)$ is gaugeable. Then for
There exist $c>1$ such that
\[
c^{-1}\le\inf_{x, y\in D} E^x_y \left[ e_q(\zeta )\right] 
\le \sup_{x, y\in D} E^x_y \left[ e_q(\zeta )\right]  \le c.
\]
\end{thm}

\begin{thm}\label{4a}
Assume that $(D, q)$ is gaugeable. 
There exists a constant  $c>1$ such that
\[
c^{-1}\, G_D(x, y)\le V_q(x, y)
\le c \, G_D(x, y) \quad \hbox{ for } x, y \in D.
\]
That is, the Green function of $D$ with respect to $L^D+q$ is comparable to
the Green function of $D$ with respect to $L^D$.
\end{thm}

\begin{thm}\label{g3:rep}
Assume that $(D, q)$ is gaugeable. 
For all $x\in D$ and $z\in\overline{D}^c$,
\[
E^x_z \left[ e_q(\zeta )\right] =
1+K_D(x, z)^{-1}\int_DV_q(x, w)q(w)K_D(w, z)dw.
\]
\end{thm}

\pf For any Borel measurable subset $A\subset  D$, by \ref{eqn:pois}
\begin{eqnarray*}
&&\int_A\frac{G_D(x, y)|q(y)|K_D(y, z)}{K_D(x, z)}dy\\
&&=
\frac{A(n, \alpha)}{K_D(x, z)}\int_AG_D(x, y)|q(y)|\left(\int_D\frac{G_D(y, w)}
{|w-z|^{n+\alpha}}dw\right)dy\\
&&=
\frac{A(n, \alpha)}{K_D(x, z)}\int_D\frac{G_D(x, w)}
{|w-z|^{n+\alpha}}\left(\int_A\frac{G_D(x, y)|q(y)|G_D(y, w)}{G_D(x, w)}dy
\right)dw.
\end{eqnarray*}
The family of functions
\[
\left\{\frac{G_D(x, \cdot)|q(\cdot)|G_D(\cdot, w)}{G_D(x, w)},
\,  x, w\in D
\right\}
\]
is uniformly integrable by (\ref{eqn:3g5}) and
 therefore the family of functions
\[
\left\{\frac{G_D(x, \cdot)|q(\cdot)|K_D(\cdot, z)}{K_D(x, z)}, 
\, x\in D,
z\in\overline{D}^c\right\}
\]
is uniformly integrable. Using Fubini's theorem, 
\begin{eqnarray*}
&&E^x_z\left[ \int^{\zeta}_0e_q(t)|q(X_t)|dt\right]\\
&=&
\int^{\infty}_0E^x_z\left[ e_q(t) |q(X_t)| ; \, t<\zeta \right]dt\\
&=& K_D(x, z)^{-1}\int^{\infty}_0 E^x\left[ e_q(t)|q(X_t)|K_D(X_t, z); \,
t<\tau_D
\right]dt\\
&=& K_D(x, z)^{-1}\int^{\infty}_0\int_Du_q(t, x, y)|q(w)|K_D(w, z)dwdt\\
&=& K_D(x, z)^{-1}\int_DV_q(x, w)|q(w)|K_D(w, z)dw <\infty .
\end{eqnarray*}
Hence 
\begin{eqnarray*}
E^x_z[e_q(\zeta)]&=&  1+E^x_z\left[ \int^{\zeta}_0e_q(t)
q(X_t)dt\right] \\
&=&1+K_D(x, z)^{-1}\int_DV_q(x, w)q(w)K_D(w, z)dw.
\end{eqnarray*}
\qed

\begin{thm}\label{4b}
Assume that $(D, q)$ is gaugeable.
There exists a constant  $c>1$ such that
\[
c^{-1}\le\inf_{(x, z)\in D\times {\overline D}^c}E^x_z \left[e_q(\tau_D)\right]
\le\sup_{(x, z)\in D\times {\overline D}^c}E^x_z \left[e_q(\tau_D) \right]
\le c.
\]
\end{thm}

\pf
By (\ref{eqn:pois})
\begin{eqnarray*}
E^x_z \left[ e_q(\tau_D) \right] &=& 1+\frac1{K_D(x, z)}\int_DV_q(x, w)
q(w)\int_D
\frac{G_D(w, v)}{|v-z|^{n+\alpha}}dvdw\\
&&=\frac{A(n, \alpha)}{K_D(x, z)}\int_D\frac{G_D(x, v)}{|v-z|^{n+\alpha}}
\left(1+ \frac1{G_D(x, v)}\int_DV_q(x, w)q(w)G_D(w, v)dw \right) dv\\
&&=\frac{A(n,\alpha)}{K_D(x, z)}\int_D\frac{G_D(x, v)}{|v-z|^{n+\alpha}}
   E^x_v \left[ e_q(\zeta )\right]  dv.
\end{eqnarray*}
The assertion of this theorem then follows from Theorem \ref{co:gau}.
\qed

\begin{thm}\label{4c}
Assume that $(D, q)$ is gaugeable.
There exists a constant  $c>1$ such that
\[
c^{-1}\le\inf_{(x, z)\in D\times \partial D}E^x_z \left[e_q(\tau_D)\right]
\le\sup_{(x, z)\in D\times \partial D}E^x_z \left[e_q(\tau_D) \right]
\le c.
\]
\end{thm}

\pf It follows from the (\ref{eqn:3g6}) that the family of functions
\[
\left\{\frac{G_D(x, \cdot)|q(\cdot)|M_D(\cdot, z)}{M_D(x, z)}, x\in D,
z\in\partial D\right\}
\]
is uniformly integrable. A similar calculation as that given in the
proof of Theorem \ref{g3:rep} yields
\[
E^x_z[e_q(\zeta)]=1+\frac1{M_D(x, z)}\int_DV_q(x, u)q(u)M_D(u, z)du.
\]
Let $x_0$ be a fixed point in $D$. Then from \cite{s3:che}
we know that
\begin{eqnarray*}
\lim_{y\to z}\frac1{G_D(x, y)}V_q(x, w)q(w)G_D(w, y)
&=&\lim_{y\to z}V_q(x, w)q(w) \frac{G_D(w, y)/
G_D(x_0, y)}
{G_D(x, y)/G_D(x_0, y)}\\
&=&V_q(x, w)q(w) \frac{M_D(w, z)}{M_D(x, z)}.
\end{eqnarray*}
Now using the uniform integrability of the family of functions
\[
\left\{\frac{G_D(x, \cdot)|q(\cdot)|G_D(\cdot, y)}{G_D(x, y)},
\,  x, y\in D
\right\}
\]
and letting $y\to z\in\partial D$ in Theorem \ref{ga:rep} yield
\[
\lim_{y\to z}E^x_y[e_q(\zeta)]=
1+\frac1{M_D(x, z)}\int_DV_q(x, u)q(u)M_D(u, z)du.
\]
Therefore
\begin{equation}
\lim_{y\to z}E_x^y[e_q(\zeta)]=
E_x^z[e_q(\zeta)]
\end{equation}
and the theorem now follows from Theorem \ref{co:gau}.
\qed

Combining Theorems \ref{co:gau}, \ref{4b} and \ref{4c} we get the following

\begin{thm} (Conditional Gauge Theorem).
Assume that $(D, q)$ is gaugeable. 
There exists a constant  $c>1$ such that
\[
c^{-1}\le\inf_{(x, z)\in D\times \RR^n}E^x_z \left[e_q(\zeta )\right]
\le\sup_{(x, z)\in D\times \RR^n}E^x_z \left[e_q(\zeta ) \right]
\le c.
\]
\end{thm}

For $q\in \KK_{n, \alpha}$, let $\psi_0>$ be the ground
state of $L^D+q$ and let $u_q(t, x, y)$ be the density kernel
of $L^D+q$ with respect to the Lebesgue measure in $D$.

\begin{thm}\label{thm:eigen}
Suppose that $(D, q)$ is gaugeable.
\begin{description}
\item{(1)} There is a constant $c>1)$ such that
$$ c^{-1}\phi_0 (x) \leq \psi_0(x) \leq c \phi_0(x) \quad 
  \hbox{for all } x\in D;
$$
\item{(2)} For each $t>0$ there is a constant $c_t>1$ such that
$$c_t^{-1} p^D(t, x, y) \leq u_q(t, x, y) \leq c_t\, p^D(t, x, y)
\quad \hbox{for all } x, y \in D .
$$
\end{description}
\end{thm}

\pf The proof is the same as  that for Theorem 2 in
Ba\~nuelos \cite{Ba} and is thus omitted.   \qed

\medskip

\noindent {\bf Remark.} It is possible to extend conditional gauge theorem beyond bounded Lipschitz domains to domains having the following
properties. Suppose that $D$ is a domain having finite
Lebesgue measure and $q\in \KK_{n, \alpha}$
such that $L^D+q$ is intrinsic ultracontractive.
Assume also that $q$ admits a decomposition $q=q_1+q_2$
with $q_1$ and $q_2$ satisfying (\ref{eqn:kato}). Then 
Theorem \ref{g1:cmp} remains true by exactly the same argument.
Assume further that
\begin{equation}\label{eqn:unif}
\left\{{G_D (x, \cdot)|q(\cdot)|G_D(\cdot, y)\over G_D(x, y)}
  : \ x, y\in D\right\}
\end{equation}
is uniformly integrable on $D$. Then Theorems \ref{gn:rel}--\ref{4a}
hold by the same proof as those for  Theorems 5.3--5.7
in \cite{s2:che}. In this case, Theorem \ref{thm:eigen} holds as well.
When $D$ is a bounded Lipschitz domain, the above  conditions are satisfied
due to 3G Theorem \ref{2b}.

\vspace{.5in}
\begin{singlespace}

\end{singlespace}
\end{doublespace}
\end{document}